\documentclass{birkau}

\usepackage{amssymb, latexsym, url, mathtools}

\newcommand{\Om}{\Omega}
\newcommand{\mb}{\mathbb}
\newcommand{\mbN}{\mathbb{N}}

\newcommand{\mc}{\mathcal}
\newcommand{\mcV}{\mathcal{V}}
\newcommand{\mcW}{\mathcal{W}}
\newcommand{\mcS}{\mathcal{S}}
\newcommand{\mcL}{\mathcal{L}}
\newcommand{\mcU}{\mathcal{U}}

\newcommand{\wt}{\widetilde}

\newcommand{\mbf}{\mathbf}
\newcommand{\mf}{\mathsf}
\newcommand{\bfx}{\mathbf{x}}
\newcommand{\bfy}{\mathbf{y}}
\newcommand{\bfz}{\mathbf{z}}
\newcommand{\bfa}{\mathbf{a}}
\newcommand{\bfb}{\mathbf{b}}
\newcommand{\bfc}{\mathbf{c}}
\newcommand{\bfd}{\mathbf{d}}

\theoremstyle{plain}
\newtheorem{theorem}{Theorem}[section]
\newtheorem{corollary}[theorem]{Corollary}
\newtheorem{proposition}[theorem]{Proposition}
\newtheorem{lemma}[theorem]{Lemma}

\theoremstyle{definition}
\newtheorem{definition}[theorem]{Definition}

\newtheorem{example}[theorem]{Example}

\numberwithin{equation}{section}

\begin{document}

\title[Mal'tsev products of varieties]{Mal'tsev products of varieties, II}

\author[T. Penza]{Tomasz Penza}
\address{Faculty of Mathematics and Information Science\\
Warsaw University of Technology\\
Warsaw, Poland}
\email{T.Penza@mini.pw.edu.pl}

\author[A.B. Romanowska]{Anna B. Romanowska}
\address{Faculty of Mathematics and Information Science\\
Warsaw University of Technology\\
Warsaw, Poland}
\email{A.Romanowska@mini.pw.edu.pl}



\begin{abstract}
The Mal'tsev product of two varieties of the same similarity type is not in general a variety, because it can fail to be closed under homomorphic images. In the previous paper we provided a new sufficient condition for such a product to be a variety. In this paper we extend that result by weakening the assumptions regarding the two varieties. We also explore the various special cases of our new result and provide a number of examples of its application.
\end{abstract}

\maketitle

\section{Introduction}

This paper is a continuation of the paper~\cite{PR21} by the same authors. The reader should consult~\cite{PR21} for further background and all notions that are not explicitly defined.

Let $\mcV$ and $\mcW$ be varieties of the same similarity type $\tau \colon \Omega \rightarrow \mathbb{N}$. The  \emph{Mal'tsev product} $\mc{V} \circ \mc{W}$ of $\mc{V}$ and $\mc{W}$ consists of all algebras $A$ of type $\tau$ with a congruence $\theta$, such that $A/\theta$ belongs to $\mc{W}$ and every congruence class of~$\theta$ that is a subalgebra of $A$ belongs to $\mc{V}$. Each algebra $A$ of type~$\tau$ has the smallest congruence such that the corresponding quotient algebra belongs to~$\mcW$. This congruence is called its $\mcW$-\emph{replica congruence} and will be denoted by~$\varrho$. (See e.g.~\cite[Ch.~3]{RS02}.) The congruence $\theta$ in the definition of the Mal'tsev product may be taken to be the $\mcW$-replica congruence $\varrho$ of $A$. (See~\cite{M67}.) Thus the definition of the Mal'tsev product of varieties becomes
\begin{equation}
\mc{V} \circ \mc{W} = \{A \mid (\forall{a \in A}) \ (a/\varrho \leq A \, \Rightarrow \, a/\varrho \in \mcV)\}.
\end{equation}
By results of Mal'tsev~\cite[Ths.~1, 2]{M67}, it is known that the Mal'tsev product $\mc{V} \circ \mc{W}$ is closed under the formation of subalgebras and of direct products. However in general, it is not closed under homomorphic images. We are interested in sufficient conditions for the Mal'tsev product $\mc{V} \circ \mc{W}$ to be a variety.

If the factor $\mcW$ of the Mal'tsev product $\mcV \circ \mcW$ is idempotent, then each $\mcW$-replica congruence class $a/\varrho$ of any algebra $A$ in $\mcV \circ \mcW$ is a subalgebra of $A$, and
\[
\mcV \circ \mcW = \{A \mid (\forall{a \in A})\,\, (a/\varrho \in \mcV)\}.
\]
In this case, we say that $A$ is a \emph{$\mcW$-sum of $\mcV$-algebras}. (See~\cite[Sec. 1]{PR21}). In this paper we extend the main result of~\cite{PR21}, in which the second factor~$\mcW$ of the Mal'tsev product $\mcV \circ \mcW$ is required to be an idempotent variety. Now it is  allowed to be a member of a wider class of varieties that we named \emph{term idempotent varieties}. This change forces us to pay more attention to the role of idempotent elements (or just \emph{idempotents}) in the theory of Mal'tsev products. Recall that for an algebra $A$ with a congruence $\theta$, a congruence class $a/\theta$ is a subalgebra of $A$ iff $a/\theta$ is an idempotent element of the quotient algebra $A/\theta$. Thus for $A \in \mcV \circ \mcW$, the congruence classes $a/\varrho$ of the $\mcW$-replica congruence $\varrho$ that are subalgebras of $A$ are precisely those that are idempotents of $A/\varrho \in \mcW$.
Some special terms are used to keep track of those congruence classes.
Let $\mcV$ be a variety of type $\tau$ and $t$ be a term of this type. If $\mcV$ satisfies the identities
\begin{equation}\label{E:tid}
\omega(t,\dots,t) = t,
\end{equation}
for all basic operations $\omega \in \Om$, then we say that $t$ is a \emph{term idempotent} of~$\mcV$. (See~\cite[Sec. 1]{PR21}.) (Iskander \cite{I84a} and \cite{I84b} uses the name \emph{unit term} for a unary term idempotent.)
To justify the name, first recall that the free $\mcV$-algebra over $X$ may be represented as the quotient $X\Om/\varrho$, where $X\Om$ is the absolutely free algebra $X\Om$ over $X$ and $\varrho$ is its $\mcV$-replica congruence. (See e.g.~\cite[Ch.~3]{RS02}.)
Then note that a term $t$, with variables in a set $X$, is a term idempotent of $\mcV$ precisely if $t/\varrho$ is an idempotent of the free $\mcV$-algebra over $X$.
A motivating example is provided by the term $t(x)\coloneqq xx^{-1}$ in a variety of groups or of inverse semigroups. Note that a variety $\mcV$ is idempotent if a variable $x$ is a term idempotent of $\mcV$. Moreover, if a variety $\mcV$ is idempotent, then every term of type $\tau$ is a term idempotent of $\mcV$. Note also that if $t$ is a term idempotent of $\mcV$ and $A \in \mcV$, then for each $a \in A$, the element $t(a)$ is an idempotent of $A$.

We will restrict our attention to types with no symbols of nullary operations. It is a reasonable assumption when dealing with Mal'tsev products for the following reasons.
First, note that if the type of $\mcV$ and $\mcW$ contains symbols of nullary operations, then each algebra $A$ in $\mcV \circ \mcW$ has only one congruence class of $\varrho$ that is a subalgebra, namely the one containing the constants. If additionally $\mcW$ is idempotent, then $\mcV \circ \mcW$ is just the variety $\mcV$. Then, if the type
contains a symbol $c$ of a nullary operation, then one can replace it by a symbol of a constant unary basic operation $c(x)$, and in this way obtain equivalent varieties $\mcV'$ and $\mcW'$ of a type without constants. The varieties $\mcV'$ and $\mcW'$ satisfy the identity $c(x) = c(y)$, and the unary operation is constant on all algebras of these varieties.

This paper is organized as follows.
Section~\ref{S:summary} contains a summary of earlier sufficient conditions for a Mal'tsev product $\mcV \circ \mcW$ to be a variety. In Section~\ref{S:termidemp} we introduce term idempotent varieties and investigate their properties. Section~\ref{S:Main} contains the main result of this paper, Theorem~\ref{T:newmain}, providing a sufficient condition for the Mal'tsev product of a variety $\mcV$ and a term idempotent variety $\mcW$ to be a variety. This theorem extends several earlier results and has a number of interesting consequences and applications that are discussed in Section~\ref{S:corandex}. Finally, in Section~\ref{S:purelypolar}, we
investigate a subclass of term idempotent varieties $\mcW$ consisting of the so-called \emph{polarized} varieties. They have a rather special property that $\mcV \circ\mcW$ is a variety for any variety~$\mcV$. This result is not a consequence of the main theorem of Section~\ref{S:Main}.

With the exception of some examples, we usually assume that $\mcV$ and $\mcW$ are varieties of the same finitary type $\tau \colon \Omega \rightarrow \mbN$ without symbols of nullary operations, and that all varieties, algebras and terms are of this type. If a variety $\mcV$ satisfies an identity $p = q$, then the terms $p$ and $q$ are called \emph{equivalent} in $\mcV$ or $\mcV$-\emph{equivalent}. An identity is \emph{trivial} if it is of the form $p = p$. A term $t(x_1,\dots,x_n)$ is called \emph{constant} in $\mcV$, if $\mcV$ satisfies the identity $t(x_1,\dots,x_n) = t(y_1,\dots,y_n)$.

We usually abbreviate lists of variables $x_1,\dots, x_n$ as $\bfx$. For a term $t = t(x_1,\dots,x_n)$ we also write $t(\bfx)$. Note that $t(\bfx)$ need not necessarily involve the full set $x_1,\dots,x_n$ of variables from $\mathbf x$. Similarly, we abbreviate lists of elements $a_1,\dots,a_n$ of some algebra as $\bfa$, and we write $t(\bfa)$ for $t(a_1,\dots ,a_n)$.
With the exception of this special notation, we follow the usage of
notation and conventions similar to those of~\cite{BPR20}, ~\cite{PR21} and~\cite{RS02}.

For further information regarding Mal'tsev products, we refer the reader to~\cite{M67} and \cite{M71}. For universal algebra, see~\cite{CB12} and~\cite{RS02}.

\section{A brief summary of earlier results}\label{S:summary}

We proceed with a brief summary of the earlier sufficient conditions for $\mcV \circ \mcW$ to be a variety. However let us first recall a result about the identities true in $\mcV \circ \mcW$.

\begin{definition}\cite[Def.~2.1]{PR21}
Let $\mcV$ and~$\mcW$ be varieties of type $\tau$, and let $\Sigma$ be an equational base for~$\mcV$. We define the following set $\Sigma^{\mcW}$ of identities:
\begin{align*}
\Sigma^{\mcW} \coloneqq \{&u(r_1, \dots, r_n) = v(r_1, \dots, r_n)  \mid\\
    & (u = v) \in \Sigma, \\
    &\forall \, i = 1,\dots ,n-1, \ \ \mcW \models r_i = r_{i+1}, \\
    &\forall \, \omega \in \Om, \ \ \mcW \models \omega(r_1, \dots, r_1) = r_1\}.
\end{align*}
\end{definition}
\noindent The last two conditions of this definition imply that $\mcW \models \omega(r_i, \dots, r_i) = r_i$ for all $\omega \in \Omega$ and each $i = 1, \dots ,n-1$, and thus all $r_i$ are term idempotents of $\mcW$.
One can say that every identity in $\Sigma^{\mcW}$ is obtained from some identity of~$\Sigma$ by substituting for its variables pairwise $\mcW$-equivalent term idempotents of $\mcW$.

\begin{theorem}\cite[Lem.~2.2]{PR21}\label{T:thevar}
Let $\mcV$ and $\mcW$ be varieties of type $\tau$, and let $\Sigma$ be an equational base for $\mcV$. Then the variety $\mf{H}(\mcV \circ \mcW)$ generated by the Mal'tsev product $\mcV \circ \mcW$ is defined by the identities $\Sigma^{\mcW}$.
\end{theorem}

\begin{corollary}
If the Mal'tsev product $\mcV \circ \mcW$ is a variety, then $\Sigma^{\mcW}$ is an equational base for $\mcV \circ \mcW$.
\end{corollary}

\begin{theorem}\cite[Th.~3.3]{PR21}\label{T:main}
Let $\mcV$ and $\mcW$ be varieties of type $\tau$, and let $\mcW$ be idempotent.
If there exist terms $f(x,y)$ and $g(x,y)$ such that
\begin{itemize}
\item[(a)] $\mcV \models f(x,y) = x\ \ \mbox{and}\ \ \mcV \models g(x,y) = y$,
\item[(b)] $\mcW \models f(x,y) = g(x,y)$,
\end{itemize}
then the Mal'tsev product $\mcV \circ \mcW$ is a variety.
\end{theorem}

An identity is called \emph{regular} if both its sides contain precisely the same variables; otherwise it is \emph{irregular}. Furthermore an identity is called \emph{strongly irregular}, if it is of the form $t(x,y) = x$, where $t(x,y)$ is a binary term containing both variables $x$ and $y$. A variety is called \emph{strongly irregular} if it satisfies a strongly irregular identity. For a \emph{plural} type $\tau$, i.e. one
with no nullary operations and at least one non-unary operation, the variety $\mcS_{\tau}$ of $\tau$-semilattices is the unique variety of type $\tau$ that is equivalent to the variety $\mcS$ of semilattices. This variety satisfies precisely all the regular identities of type $\tau$. (See~\cite{PR92} and~\cite{BPR20} for details.)
As a corollary of Theorem~\ref{T:main} one obtains the following theorem.

\begin{theorem}\cite[Th.~6.3]{BPR20}\label{T:BPR}
If $\mcV$ is a strongly irregular variety of a plural type~$\tau$, then $\mcV \circ \mcS_{\tau}$ is a variety.
\end{theorem}

\noindent Algebras in $\mcV \circ \mcS_\tau$ are called \emph{semilattice sums of $\mcV$-algebras}.

The main result of this paper is a common generalization of Theorem~\ref{T:main} and of the following theorem of Bergman.

\begin{theorem}\cite[Cor.~2.3]{B20}\label{T:BM}
If $\mcV$ and $\mcW$ are idempotent subvarieties of a congruence permutable variety, then $\mcV \circ \mcW$ is a variety.
\end{theorem}

\section{Term idempotent varieties} \label{S:termidemp}

We start this section with a special property of term idempotents. The set $X\Om$ of terms of a given type $\tau$ (without constants) over a countably infinite set $X$ of variables is preordered by the following relation: $p(x_1 \dots x_n) \preceq q$ iff there exist terms $t_1, \dots ,t_n$ of type $\tau$ such that $q = p(t_1 \dots t_n)$. Note that if $p \preceq q$ and $q \preceq p$, then $p$ and $q$ are the same to within a renaming of the variables. Using the rules of equational logic one easily obtains the following lemma.

\begin{lemma}\label{L:preorder}
If $p$ is a term idempotent of a variety $\mcV$, and $p \preceq q$, then $q$ is also a term idempotent of $\mcV$.
\end{lemma}

It is known that for a preordered set $(P,\preceq)$, the relation $\alpha$ defined on $P$ by
\[
(p,q) \in \alpha \ \mbox{iff} \ p \preceq q \ \mbox{and} \ q \preceq p
\]
is an equivalence relation. Furthermore, the relation $\leq$ defined on $P/\alpha$ by
\[
p/\alpha \leq q/\alpha \ \mbox{iff} \ p \preceq q
\]
is an order relation.
An upper set of a preordered set $(P,\preceq)$ can be defined similarly as in the case of an ordered set. A subset $Q$ of $P$ is an \emph{upper set}, if whenever $p \in Q$, $q \in P$ and $p \preceq q$, then $q \in Q$. Thus Lemma~\ref{L:preorder} shows that term idempotents of a given variety form an upper set of $(X\Om,\preceq)$.
It is easy to see that
the variables of~$X$ form one class of $\alpha$. This class is the minimum of the ordered set $(X\Om/\alpha,\leq)$ and obviously each variable of~$X$ is related by $\preceq$ with any other element of $X\Om$.
Thus if a variable is a term idempotent, then all terms are term idempotents. In other words, a variety is idempotent iff the upper set of its term idempotents contains all terms.

\begin{definition}\label{D:termidemvar}
A nontrivial identity $p = q$ satisfied in a variety $\mcV$ will be called \emph{term idempotent}, if both $p$ and $q$ are term idempotents of $\mcV$. A variety $\mcV$ will be called \emph{term idempotent}, if every nontrivial identitity it satisfies is term idempotent.
\end{definition}

Note that every idempotent variety is term idempotent.
Below we provide some examples of term idempotent varieties that are not idempotent.

\begin{example} \label{Ex:cs}
Let $\mc{CS}$ be the variety of constant semigroups, i.e. the variety of groupoids (magmas, binars) defined by the identity $x y = z t$. A nontrivial identity $p = q$ is satisfied in~$\mc{CS}$ precisely if neither $p$ nor $q$ is a variable. (Cf~\cite[Ex.~2.5]{PR21}.) In particular $\mc{CS}$ satisfies $p \cdot p = p$ for every term $p$ different from a variable. Consequently, all such terms are term idempotents of $\mc{CS}$, and so $\mc{CS}$ is a term idempotent variety.
\end{example}

\begin{example}\label{Ex:ca}
The variety $\mc{C}_\tau$ of \emph{constant algebras} of type $\tau$ is defined by the identities
\[
\omega(x_1,\dots,x_n) = \varphi(y_1,\dots,y_m),
\]
for all $\omega,\varphi \in \Om$. It satisfies precisely the nontrivial identities whose neither side is a variable. Every term different from a variable is constant in $\mc{C}_\tau$. If type $\tau$ consists of a single binary operation, then $\mc{C}_\tau$ is just the variety $\mc{CS}$. For any other type, $\mc{C}_\tau$ is equivalent to $\mc{CS}$. An argument analogous to that of Example~\ref{Ex:cs} shows that $\mc{C}_\tau$ is a term idempotent variety.
\end{example}

\begin{example} \label{Ex:RS}
Let $\mc{RS}$ be the variety of semigroups defined by the identities
\begin{equation} \label{E:rs}
        (xy)z = xz = x(yz).
\end{equation}
The subvariety of $\mc{RS}$ defined by the idempotent law $xx = x$ is the variety of \emph{rectangular bands}. Algebras in $\mc{RS}$ will be called \emph{rectangular semigroups}. Recall that in any variety of semigroups each term $t$ is equivalent to a product of variables $t = x_1\cdots x_n$. If $t$ is different from a variable, then~\eqref{E:rs} implies that $\mc{RS}\models t = x_1 x_n$. Thus $\mc{RS}$ satisfies
\[
    t \cdot t = x_1 x_n x_1 x_n = x_1 x_n = t.
\]
Consequently, all terms different from a variable are term idempotents of $\mc{RS}$. Since all nontrivial identities derivable from~\eqref{E:rs} have both sides different from a variable, it follows that $\mc{RS}$ is term idempotent.
\end{example}

\begin{example}\label{Ex:Un}
For $n \geq 0$, let $\mcU_n$ be the variety of monounary algebras $(A,f)$ defined by the identity
\begin{equation} \label{E:Un}
        f(f^n(x))=f^n(x).
\end{equation}
Clearly, $f^n(x)$ is a term idempotent of $\mcU_n$. Recall that each term of monounary type has the form $f^m(x)$ for some $m \geq 0$ (with $f^0(x)$ being just $x$). Then by~\eqref{E:Un}, if $m \geq n$, then $\mcU_n \models f^m(x) = f^n(x)$. Every nontrivial identity derivable from~\eqref{E:Un} is of the form $f^k(x) = f^l(x)$ for different $k,l \geq n$. So both sides of such an identity are term idempotents of $\mcU_n$, and hence $\mcU_n$ is term idempotent. Observe that $\mcU_0$ is idempotent, in $\mcU_1$ all terms different from a variable are term idempotent, and if $n \geq 2$, then in $\mcU_n$ not all terms different from a variable are term idempotents.
\end{example}

Some regular varieties provide further examples of term idempotent varieties. Recall that the \emph{regularization} $\widetilde{\mcV}$ of a variety $\mcV$ of a plural type $\tau$ is the variety defined by all the regular identities satisfied in $\mcV$. Equivalently, $\widetilde{\mcV}$  can be defined as the join $\mcV \vee \mcS_\tau$ of $\mcV$ and the variety $\mcS_\tau$ of $\tau$-semilattices. It is known that if $\mcV$ is irregular, then each algebra in $\widetilde{\mcV}$ is a semilattice sum of $\mcV$-algebras. If $\mcV$ is strongly irregular, then $\widetilde{\mcV}$ coincides with the class of Płonka sums of $\mcV$-algebras. (See e.g.~\cite{P69},~\cite[Ch.~4]{RS02},~\cite{PR92}.)

\begin{proposition} \label{P:regularization}
Let $\mcV$ be a variety of a plural type $\tau$. If $\mcV$ is term idempotent, then $\widetilde{\mcV}$ is also term idempotent.
\end{proposition}
\begin{proof}
Let $u = v$ be a nontrivial identity satisfied in $\widetilde{\mcV}$. Then $u = v$ is also satisfied in $\mcV$, and hence $u$ and $v$ are term idempotents of $\mcV$. Thus $\mcV$ satisfies the identities $\omega(u,\dots ,u) = u$ and $\omega(v,\dots ,v) = v$ for all $\omega \in \Om$. Since these identities are regular, they are also satisfied in $\widetilde{\mcV}$. Therefore $u$ and $v$ are term idempotents of $\widetilde{\mcV}$, and so  $\widetilde{\mcV}$ is a term idempotent variety.
\end{proof}

Note that none of the examples of term idempotent varieties which are not idempotent that we provided so far are strongly irregular. The next proposition shows that this is not a coincidence.

\begin{proposition}
If a variety $\mcV$ is term idempotent and strongly irregular, then $\mcV$ is idempotent.
\end{proposition}
\begin{proof}
The variety $\mcV$ satisfies a strongly irregular identity $t(x,y)= x$. This identity is nontrivial, so it is term idempotent, and thus in particular its right-hand side $x$ is a term idempotent of $\mcV$. Therefore $\mcV$ is idempotent.
\end{proof}

We conclude this section with a characterization of term idempotent varieties in terms of replica congruences. First recall a very useful description of a replica congruence which will also be used in the proof of the main result in Section~\ref{S:Main}.

\begin{definition}\label{D:ro0}
Let $\mcW$ be a variety, and let $A$ be an algebra of the same type as $\mcW$. We define a binary relation $\varrho^0$ on the universe of $A$ as follows:
$(a, b) \in \varrho^0$ if and only if there are an identity $p(\bfx) = q(\bfx)$ satisfied in $\mcW$, and elements $\bfd$ of $A$, such that $a = p(\bfd)$ and $b = q(\bfd)$.
\end{definition}

\noindent Note that the relation $\varrho^0$ is reflexive and symmetric.

\begin{proposition}\cite[Prop.~3.2]{PR21} \label{L:repcong}
Let $\mcW$ be a variety, and let $A$ be an algebra of the same type as $\mcW$. The $\mcW$-replica congruence relation $\varrho$ of $A$ coincides with the transitive closure of $\varrho^0$.
\end{proposition}

\begin{proposition} \label{P:singleelement}
Let $\mcW$ be a variety of type $\tau$. Then $\mcW$ is term idempotent if and only if, for every algebra $A$ of type $\tau$, every congruence class $a/\varrho$ of the $\mcW$-replica congruence of $A$ which is not an idempotent of $A/\varrho$, is a singleton.
\end{proposition}
\begin{proof}
$(\Rightarrow)$ Assume that $\mcW$ is a term idempotent variety. Let $a/\varrho$ be a congruence class with more than one element. We will show that $a/\varrho$ is an idempotent of $A/\varrho$. Let $b \in a/\varrho$ be an element different from $a$. By Proposition~\ref{L:repcong}, $\varrho$ is the transitive closure of $\varrho^0$. Since $(a,b) \in \varrho$, there is an element $c \in a/\varrho$ different from $a$, such that $(a,c) \in \varrho^0$. This means that there is a nontrivial identity $p(\bfx) = q(\bfx)$ true in $\mcW$, and elements $\bfd$ of $A$, such that $a = p(\bfd)$ and $c = q(\bfd)$.
Since $\mcW$ is a term idempotent variety, $p$ and $q$ are term idempotents of $\mcW$. Thus $a$, being a value of a term idempotent, is an idempotent of $A$. It follows that for each $\omega \in \Om$,
\[
a = \omega(a, \dots, a),
\]
and hence
\[
a/\varrho = \omega(a/\varrho, \dots, a/\varrho).
\]
Therefore $a/\varrho$ is an idempotent of $A/\varrho$.

$(\Leftarrow)$
Now assume that for every algebra $A$ of type $\tau$, any congruence class $a/\varrho$ which is not an idempotent of $A/\varrho$, has exactly one element.
In particular, this is true for the absolutely free algebra $X\Om$ of type $\tau$ over a countably infinite set $X$. Recall that the quotient $X\Om/\varrho$ of $X\Om$ by its $\mcW$-replica congruence $\varrho$, is the free $\mcW$-algebra over $X$. (See e.g.~\cite[Ch.~3]{RS02}.) By the definition of term idempotents, for a term $t$ of type $\tau$, the congruence class $t/\varrho$ is an idempotent of $X\Om/\varrho$ precisely if $t$ is a term idempotent of~$\mcW$. Now let $p = q$ be a nontrivial identity satisfied in $\mcW$. Then $p$ and $q$ are different elements of the same congruence class $C \coloneqq p/\varrho = q/\varrho$. Since $C$ is not a singleton, it is an idempotent of $X\Om/\varrho$. Hence $p$ and $q$ are term idempotents of $\mcW$, and therefore $\mcW$ is term idempotent.
\end{proof}

As a corollary, one obtains a result on the structure of algebras in $\mcV \circ \mcW$ for a term idempotent variety $\mcW$.

\begin{corollary}\label{C:onesubalg}
Let $\mcV$ and $\mcW$ be varieties of the same type, and let $\mcW$ be term idempotent. If $A \in \mcV \circ \mcW$, then each congruence class of the $\mcW$-replica congruence $\varrho$ of $A$ is either a subalgebra of $A$ or a singleton.
\end{corollary}

In the definition of a term idempotent variety $\mcV$ we require that all nontrivial identities true in $\mcV$ are term idempotent. One might wonder if this property is equivalent to the requirement that the set of identities used to define $\mcV$ be term idempotent. This is not the case however, since a term idempotent identity may entail nontrivial identities that are not term idempotent. As an example consider the identity $xx^{-1} = yy^{-1}$ true in the variety of groups. Both of its sides are term idempotents. However this identity implies the nontrivial identity $(xx^{-1})z = (yy^{-1})z$ whose both sides are equivalent to~$z$ which is not a term idempotent.

\section{A new sufficient condition for $\mcV \circ \mcW$ to be a variety}\label{S:Main}

We are now ready to state and prove our generalization of Theorem~\ref{T:main}.

\begin{theorem}\label{T:newmain}
Let $\mcV$ and $\mcW$ be varieties of the same type, and let $\mcW$ be term idempotent.
If there exist terms $f(x,y,z)$ and $g(x,y,z)$ such that
\begin{itemize}
\item[(a)] $\mcV \models f(x,y,y) = x\ \mbox{and}\ \mcV \models g(x,x,y) = y$,
\item[(b)] $\mcW \models f(x,x,y) = g(x,x,y)$,
\item[(c)] $f(x,x,y)$ is a term idempotent of $\mcW$,
\end{itemize}
then the Mal'tsev product $\mcV \circ \mcW$ is a variety.
\end{theorem}
\begin{proof}
We need to show that $\mf{H}(\mcV \circ \mcW) \subseteq \mcV \circ \mcW$. The proof will be divided into several parts. In what follows we assume that $A \in \mf{H}(\mcV \circ \mcW)$, i.e. $A$ is a quotient of an algebra belonging to $\mcV \circ \mcW$.

\smallskip
\textbf{A}. \emph{The $\mcW$-replica congruence $\varrho$ of $A$ coincides with the relation $\varrho^0$.}

\smallskip
By Proposition~\ref{L:repcong}, we have to show that the relation $\varrho^0$ is transitive.
Let $a, b, c \in A$ and $a\ \varrho^0\ b\ \varrho^0\ c$. If either two of the elements $a, b, c$ are equal, then by reflexivity and symmetry of $\varrho^0$ one obtains $a\ \varrho^0\ c$. So let us assume that these elements are pairwise different. Then there exist nontrivial identities $p_1(\bfx_1) = q_1(\bfx_1)$ and ${p_2(\bfx_2) = q_2(\bfx_2)}$ satisfied in $\mcW$, and sets of elements $\bfd_1$ and $\bfd_2$ in $A$, such that
\begin{align*}
        a = p_1(\bfd_1),\ \ &b = q_1(\bfd_1),\\
       &b = p_2(\bfd_2),\ \ \ c = q_2(\bfd_2).
\end{align*}
Since $\mcW$ is term idempotent, terms $p_1, q_1, p_2, q_2$ are term idempotents of $\mcW$. Let
\begin{equation*}
    p(\bfx_1,\bfx_2) \coloneqq f(p_1,q_1,p_2)\ \ \ \text{and}\ \ \ q(\bfx_1,\bfx_2) \coloneqq g(q_1,q_1,q_2).
\end{equation*}
By (b), $\mcW$ satisfies the identity $p = q$. By Theorem \ref{T:thevar}, $A$ satisfies the identities $f(p_1,q_1,q_1) = p_1$ and $g(p_2,p_2,q_2) = q_2$. It follows that
\begin{align*}
    p(\bfd_1,\bfd_2) & = f(p_1(\bfd_1),q_1(\bfd_1),p_2(\bfd_2))\\
    & = f(p_1(\bfd_1),q_1(\bfd_1),q_1(\bfd_1)) = p_1(\bfd_1) = a,
\end{align*}
and similarly,
\begin{align*}
    q(\bfd_1,\bfd_2) & = g(q_1(\bfd_1),q_1(\bfd_1),q_2(\bfd_2))\\
    & = g(p_2(\bfd_2),p_2(\bfd_2),q_2(\bfd_2)) = q_2(\bfd_2) = c.
\end{align*}
Thus $a\ \varrho^0\ c$, and hence $\varrho^0$ is transitive.

\smallskip
\textbf{B}. \emph{If $C$ is a congruence class of $\varrho$ which is a subalgebra of $A$, then $C$ satisfies the identities of (a).}

\smallskip
Let $u(x,y) = v(x,y)$ be an identity satisfied in $\mcV$. If~$C$ has only one element $a$, then $u(a,a) = a = v(a,a)$. So $C$ satisfies the identity~$u = v$. Now assume that $C$ has more than one element. Let $a,b\in C$ with $a \neq b$.
Then $(a,b)\in\varrho^0$, so there exist a nontrivial identity $p(\bfx) = q(\bfx)$ satisfied in $\mcW$ and elements~$\bfd$ of $A$, such that $a = p(\bfd)$ and $b = q(\bfd)$. The terms $p$ and $q$ are term idempotents of $\mcW$. Hence by Theorem~\ref{T:thevar}, $A$ satisfies the identities $u(p,q) = v(p,q)$ and $u(p,p) = v(p,p)$. Therefore
\begin{align*}
    &u(a,b) = u(p(\bfd),q(\bfd)) = v(p(\bfd),q(\bfd)) = v(a,b),\\
    &u(a,a) = u(p(\bfd),p(\bfd)) = v(p(\bfd),p(\bfd)) = v(a,a).
\end{align*}
It follows that $C$ satisfies any identity in at most two variables valid in $\mcV$. In particular
\begin{equation} \label{E:Cidentities}
        C \models f(x,y,y) = x\ \ \mbox{and}\ \ C \models g(x,x,y) = y.
\end{equation}

\smallskip
\textbf{C}. \emph{Assume that $\mcW \models p_i(\bfz_i) = q_i(\bfz_i)$ for $i = 1, \dots ,n-1$. Then for each $i = 1, \dots ,n$, define terms $t_{i,j}$ recursively for $j = 0, \dots ,n-1$ by}
\begin{equation} \label{E:deftij}
    t_{i,j} \coloneqq
    \begin{cases}
    p_1                      &\text{for $j = 0$,}\\
    f(q_j,p_j,t_{i,j-1})     &\text{for $0 < j < i$,}\\
    g(q_j,q_j,t_{i,j-1})     &\text{for $j \geq i$.}
    \end{cases}
\end{equation}
\emph{Set
\begin{equation} \label{E:defti}
    t_i \coloneqq t_{i,n-1}.
\end{equation}
Then $\mcW \models t_{i} = t_{i+1}$ for $i = 1, \dots ,n-1$.}

\smallskip
Let $1 \leq i \leq n-1$. Then by~\eqref{E:deftij} we have the following equalities
\begin{align*}
&t_{i,0} = p_1 = t_{i+1,0}, \\
&t_{i,1} = f(q_1,p_1,t_{i,0}) = f(q_1,p_1,t_{i+1,0}) = t_{i+1,1}, \\
&t_{i,2} = f(q_2,p_2,t_{i,1}) = f(q_2,p_2,t_{i+1,1}) = t_{i+1,2}, \\
  & \ \, \vdots \\
&t_{i,i-1} = f(q_{i-1},p_{i-1},t_{i,i-2}) = f(q_{i-1},p_{i-1},t_{i+1,i-2}) = t_{i+1,i-1}.
\end{align*}
Since the identities $p_i = q_i$ and $f(x,x,y) = g(x,x,y)$ are valid in $\mcW$, it follows that
\[
\mcW \models t_{i,i} = g(q_i,q_i,t_{i,i-1}) = f(q_i,p_i,t_{i+1,i-1}) = t_{i+1,i}.
\]
Hence, again by~\eqref{E:deftij}
\begin{align*}
&\mcW \models t_{i,i+1} = g(q_{i+1},q_{i+1},t_{i,i}) = g(q_{i+1},q_{i+1},t_{i+1,i}) = t_{i+1,i+1}, \\
&\mcW \models t_{i,i+2} = g(q_{i+1},q_{i+1},t_{i,i+1}) = g(q_{i+1},q_{i+1},t_{i+1,i+1}) = t_{i+1,i+2}, \\
  & \ \, \vdots \\
&\mcW \models t_{i,n-1} = g(q_{n-1},q_{n-1},t_{i,n-2}) = g(q_{n-1},q_{n-1},t_{i+1,n-2}) = t_{i+1,n-1}.
\end{align*}
Consequently $\mcW \models t_{i} = t_{i+1}$.

\smallskip
\textbf{D}. \emph{The term $t_1$ is a term idempotent of $\mcW$.}

\smallskip
By (b) and (c), $g(x,x,y)$ is a term idempotent of $\mcW$. Then we have
\[
g(x,x,y) \preceq g(q_{n-1},q_{n-1},t_{1,n-2}) = t_{1,n-1} = t_1.
\]
By Lemma~\ref{L:preorder}, it follows that $t_1$ is a term idempotent of $\mcW$.

\smallskip
\textbf{E}. \emph{Let $a_1, \dots, a_n \in C$. There exist pairwise $\mcW$-equivalent term idempotents $t_1, \dots, t_n$ of $\mcW$ and elements $\bfc$ of $A$, such that $a_i = t_i(\bfc)$ for each $i = 1, \dots, n$.}

\smallskip
Let $1 \leq i \leq n-1$. Since $a_i \, \varrho^0 \, a_{i+1}$, there is an identity $p_i(\bfz_i) = q_i(\bfz_i)$ true in $\mcW$, and elements $\bfc_i$ of $A$, such that
\[
a_i = p_i(\bfc_i) \ \ \mbox{and} \ \ a_{i+1} = q_i(\bfc_i).
\]
Denote the list $\bfc_1, \dots, \bfc_n$ by $\bfc$. Define terms $t_1,\dots,t_n$ by \eqref{E:deftij} and \eqref{E:defti}. Then \textbf{C} and \textbf{D} imply that the terms $t_1, \dots, t_n$ are $\mcW$-equivalent term idempotents of $\mcW$. By \eqref{E:Cidentities}, one obtains the following equalities:
\begin{align*}
&t_{i,0}(\bfc) = p_1(\bfc_1) = a_1, \\
&t_{i,1}(\bfc) = f(q_1(\bfc_1),p_1(\bfc_1),t_{i,0}(\bfc)) = f(a_2,a_1,a_1) = a_2, \\
&t_{i,2}(\bfc) = f(q_2(\bfc_2),p_2(\bfc_2),t_{i,1}(\bfc)) = f(a_3,a_2,a_2) = a_3, \\
  & \quad \quad \vdots \\
&t_{i,i-1}(\bfc) = f(q_{i-1}(\bfc_{i-1}),p_{i-1}(\bfc_{i-1}),t_{i,i-2}(\bfc)) = f(a_i,a_{i-1},a_{i-1}) = a_i, \\
&t_{i,i}(\bfc) = g(q_i(\bfc_i),q_i(\bfc_i),t_{i,i-1}(\bfc)) = g(a_{i+1},a_{i+1},a_i) = a_i, \\
&t_{i,i+1}(\bfc) = g(q_{i+1}(\bfc_{i+1}),q_{i+1}(\bfc_{i+1}),t_{i,i}(\bfc)) = g(a_{i+2},a_{i+2},a_i) = a_i. \\
  & \quad \quad \vdots \\
&t_{i,n-1}(\bfc) = g(q_{n-1}(\bfc_{n-1}),q_{n-1}(\bfc_{n-1}),t_{i,n-2}(\bfc)) = g(a_n,a_n,a_i) = a_i.
\end{align*}
Therefore $a_i = t_{i,n-1}(\bfc) = t_i(\bfc)$ for each $i = 1,\dots,n$.

\smallskip
\textbf{F}. \emph{The subalgebra $C$ of $A$ satisfies any identity
\[
u(x_1, \dots, x_n) = v(x_1, \dots, x_n)
\]
valid in $\mcV$.}

\smallskip
Let $a_1,\dots,a_n \in C$. Let $t_1,\dots,t_n$ be terms and $\bfc$ be elements of $A$ satisfying the condition of $\textbf{E}$. By Theorem~\ref{T:thevar}, the identity $u(t_1, \dots, t_n) = v(t_1, \dots, t_n)$ is valid in $A$. Hence
\[
u(a_1, \dots, a_n) = u(t_1(\bfc), \dots, t_n(\bfc)) = v(t_1(\bfc), \dots, t_n(\bfc)) = v(a_1, \dots, a_n).
\]

\smallskip
\textbf{G}. \emph{The Mal'tsev product $\mcV \circ \mcW$ is a variety.}

\smallskip
By \textbf{F}, we conclude that $C \in \mcV$, and consequently that $A \in \mcV \circ \mcW$. Thus the inclusion $\mf{H}(\mcV \circ \mcW) \subseteq \mcV \circ \mcW$ holds, and so $\mcV \circ \mcW$ is a variety.
\end{proof}

\section{Consequences and examples} \label{S:corandex}

Theorem~\ref{T:newmain} has a number of interesting consequences. First note that since every term in an idempotent variety is a term idempotent, one easily obtains the following corollary.

\begin{corollary} \label{C:1}
Let $\mcV$ and $\mcW$ be varieties of the same type, and let $\mcW$ be idempotent.
If there exist terms $f(x,y,z)$ and $g(x,y,z)$ such that
\begin{itemize}
\item[(a)] $\mcV \models f(x,y,y) = x\ \ \mbox{and}\ \ \mcV \models g(x,x,y) = y$,
\item[(b)] $\mcW \models f(x,x,y) = g(x,x,y)$,
\end{itemize}
then the Mal'tsev product $\mcV \circ \mcW$ is a variety.
\end{corollary}

Another special case is when the terms $f(x,y,z)$ and $g(x,y,z)$ do not depend on the middle variable.

\begin{corollary} \label{C:2}
Let $\mcV$ and $\mcW$ be nontrivial varieties, and let $\mcW$ be term idempotent.
If there exist terms $f(x,y)$ and $g(x,y)$ such that
\begin{itemize}
\item[(a)] $\mcV \models f(x,y) = x\ \ \mbox{and}\ \ \mcV \models g(x,y) = y$,
\item[(b)] $\mcW \models f(x,y) = g(x,y)$,
\end{itemize}
then the Mal'tsev product $\mcV \circ \mcW$ is a variety.
\end{corollary}
\begin{proof}
We only need to show that the condition (c) of Theorem~\ref{T:newmain} is  satisfied, i.e. $f(x,y)$ is a term idempotent of $\mcW$. First note that the terms $f(x,y)$ and $g(x,y)$ cannot coincide. Otherwise, the condition (a) would imply that $\mcV \models x = y$, contradicting the nontriviality of $\mcV$. So the identity of (b) is nontrivial, and hence $f(x,y)$ is a term idempotent of $\mcW$.
\end{proof}

\begin{example}\label{Ex:SRe}
Let $\mcV$ be the variety of groupoids defined by identities
\[
(xx)y = y = y(xx),
\]
and let $\mc{RS}$ be the variety from Example~\ref{Ex:RS}. Define
\[
f(x,y) \coloneqq x(yy) \quad \mbox{and} \quad g(x,y) \coloneqq (xx)y.
\]
It is easy to see that $\mcV$ satisfies the identities $f(x,y) = x$ and $g(x,y) = y$, and $\mc{RS}$ satisfies the identity $f(x,y) = g(x,y)$. Thus, by Corollary~\ref{C:2}, the Mal'tsev product $\mcV \circ \mc{RS}$ is a variety.
\end{example}

\begin{example}\label{Ex:strregconst}
Let $\mcV$ be a strongly irregular variety of a plural type $\tau$ that satisfies a strongly irregular identity $t(x,y) = x$. By Example~\ref{Ex:ca}, the variety $\mc{C}_\tau$ of constant algebras of type $\tau$, is a term idempotent variety. Set $f(x,y) \coloneqq t(x,y)$ and $g(x,y) \coloneqq t(y,x)$. Clearly $\mcV$ satisfies the identities $f(x,y) = x$ and $g(x,y) = y$. Since neither $f(x,y)$ nor $g(x,y)$ is a variable, it follows that $\mc{C}_\tau \models f(x,y) = g(x,y)$.  Corollary~\ref{C:2} implies that the Mal'tsev product $\mcV \circ \mc{C}_\tau$ is a variety.

Replacing the variety $\mc{C}_\tau$ by its regularization $\widetilde{\mc{C}_\tau}$ one obtains a further example.
First note that, by Proposition~\ref{P:regularization}, the regularization $\widetilde{\mc{C}_\tau}$ of $\mc{C}_\tau$ is a term idempotent variety. Keep the same terms $f(x,y)$ and $g(x,y)$ as in the previous case. Then note that
the identity $f(x,y) = g(x,y)$ is regular. Hence it is also satisfied in $\widetilde{\mc{C}_\tau}$. By Corollary~\ref{C:2} again, the Mal'tsev product $\mcV \circ \widetilde{\mc{C}_\tau}$ is a variety for any strongly irregular variety $\mcV$.
\end{example}

An additional assumption that the variety $\mcW$ of Corollary~\ref{C:2} is idempotent yields Theorem~\ref{T:main} as a special case. If we further set $g(x,y) \coloneqq y$, then the condition (a) of Corollary~\ref{C:2} reduces to only one identity and one obtains the following corollary.

\begin{corollary} \label{C:3}
Let $\mcV$ and $\mcW$ be nontrivial varieties, and let $\mcW$ be term idempotent. If there exists a term $f(x,y)$ such that
\begin{itemize}
\item[(a)] $\mcV \models f(x,y) = x$,
\item[(b)] $\mcW \models f(x,y) = y$,
\end{itemize}
then the Mal'tsev product $\mcV \circ \mcW$ is a variety.
\end{corollary}

\noindent The conditions (a) and (b) of Corollary~\ref{C:3} mean that the  varieties $\mcV$ and $\mcW$ are independent. (See e.g.~\cite[\S 3.5]{RS02}.) The join $\mcV \vee \mcW$ of independent varieties $\mcV$ and $\mcW$ consists of (algebras isomorphic to) products $A \times B$ of $A \in \mcV$ and $B \in \mcW$, and obviously
$\mcV \vee \mcW \subseteq \mcV \circ \mcW$.

\begin{example}
The variety $\mc{LZ}$ of left-zero bands (defined by the identity $xy = x$) and the variety $\mc{RZ}$ of right-zero bands (defined by the identity $xy = y$) are clearly independent. So, by Corollary~\ref{C:3}, the Mal'tsev product $\mc{LZ} \circ \mc{RZ}$ is a variety. Its subvariety $\mc{LZ} \vee \mc{RZ}$ is the variety of rectangular bands.
\end{example}

If we set $f$ and $g$ in Corollary~\ref{C:1} to be the same term, then the condition (b) is trivially satisfied, and the identities of (a) become Mal'tsev identities. So the variety $\mcV$ is congruence permutable. We thus  obtain the following extension of Theorem~\ref{T:BM}.

\begin{corollary} \label{C:Maltermidem}
Let $\mcV$ be a congruence permutable variety and $\mcW$ be an idempotent variety. Then the Mal'tsev product $\mcV \circ \mcW$ is a variety.
\end{corollary}

Typical examples of congruence permutable (or \emph{Mal'tsev}) varieties are given by varieties of groups, quasigroups, loops, rings or modules. If $\mcV$ is any such variety and $\mcW$ is an idempotent variety of the same type as $\mcV$, then the Mal'tsev product $\mcV \circ \mcW$ is a variety.

\begin{example}
Now we will consider the Mal'tsev product $\mc{G} \circ \mc{L}$ of the variety of groups and the variety of lattices. To do so we first have to describe them as varieties of the same type.

Here groups are defined as algebras $(G,\cdot,^{-1})$ with one binary and one unary operations. (See e.g.~\cite[Ex.~7.6]{BPR20}.) And the variety $\mc{G}$ of groups will be considered as the variety of algebras $(G,+,\cdot,^{-1})$, satisfying the usual identities of groups and the identity $x + y = x \cdot y$.

On the other hand, lattices will be considered as algebras of the same type as $\mc{G}$ with $x^{-1} \coloneqq x$. The variety of groups is a Mal'tsev variety, and the variety of lattices is an idempotent variety. By Corollary~\ref{C:Maltermidem}, the Mal'tsev product $\mc{G} \circ \mc{L}$ is a variety.
\end{example}

Note that in those Mal'tsev products $\mcV \circ \mcW$ considered so far which are actually varieties, the factor $\mcV$ was always strongly irregular. It is natural to ask if this strongly
irregular variety $\mcV$ could be replaced by an irregular (but not
strongly irregular) variety $\mcV$. We conclude this section with an
example of a Mal'tsev product $\mcV \circ \mcW$ that is not a variety,
where the factor $\mcV$ is irregular, but not strongly irregular, and where the factor $\mcW$ is an idempotent
(and hence term idempotent) variety.
In particular, this example shows that in Theorem~\ref{T:BPR}, which is a
corollary of Theorem~\ref{T:newmain}, the assumption of strong irregularity
cannot be replaced by irregularity alone.

\begin{example}
It will be shown that the Mal'tsev product $\mc{CS} \circ \mc{S}$ fails to be a variety. We provide an example of a groupoid that belongs to $\mc{CS} \circ \mcS$, but has a quotient that does not. Let $A$ be the groupoid defined by the following table
\begin{equation*}
    \begin{tabular}{c|ccccccc}
    $\cdot$ &$a$&$e$&$b$&$f$\\
    \hline
    $a$ & $e$ & $e$ & $b$ & $f$\\
    $e$ & $e$ & $e$ & $f$ & $f$\\
    $b$& $b$ & $f$ & $f$ & $f$\\
    $f$ & $f$ & $f$ & $f$ & $f$\\
    \end{tabular}
\end{equation*}
The groupoid $A$ is a member of $\mc{CS} \circ \mc{S}$. The semilattice replica congruence $\varrho$ of~$A$ has two congruence classes $\{a,e\}$ and $\{b,f\}$ that can be seen to be constant semigroups with constant values $e$ and $f$ respectively. The congruence $\theta$ of $A$ with congruence classes $\{a\}$, $\{b\}$ and $E \coloneqq \{e,f\}$ has the quotient $A/\theta$ given by the table
\begin{equation*}
    \begin{tabular}{c|ccccc}
   $\cdot$ &$\{a\}$&$\{b\}$&$E$\\
    \hline
    $\{a\}$  & $E$ & $\{b\}$ & $E$\\
    $\{b\}$ & $\{b\}$ & $E$ & $E$\\
    $E$  & $E$ & $E$ & $E$\\
    \end{tabular}
\end{equation*}
The semilattice replica congruence of $A/\theta$ is the all relation with one congruence class containing all elements. However $A/\theta$ is not a constant semigroup, hence $A/\theta \notin \mc{CS} \circ \mcS$. Thus $\mc{CS} \circ \mcS$ is not closed under homomorphic images, and so it is not a variety.
\end{example}

\section{Purely polarized varieties} \label{S:purelypolar}

A variety $\mcV$ that has a constant unary term idempotent is
called \emph{polarized}. (See~\cite{M67}.) Such a term is called a \emph{polar term} and is denoted by $p(x)$.
The constant value $p$ of a polar term $p(x)$ in a given algebra $A \in \mcV$ is called the \emph{pole} of this algebra.
The pole of $A$ is the unique idempotent of $A$, and a class of a congruence on $A$ is a subalgebra of $A$ iff it is the congruence class of the pole. (See~\cite[Sec.~1]{PR21}.) Since the only possible value in $A$ of a term idempotent of $\mcV$ is the pole of $A$, it follows that all term idempotents of $\mcV$ are constant and pairwise $\mcV$-equivalent.
In particular every unary term idempotent is a polar term.

Examples of polarized varieties are provided by varieties of groups with a polar term $p(x) \coloneqq xx^{-1}$, varieties of loops with a polar term $p(x) \coloneqq x/x$, and varieties of rings with a polar term $p(x)\coloneqq x - x$. One more example is given by the variety $\mc{C}_\tau$ of constant algebras where every unary term different from a variable is a polar term.

A variety will be called \emph{purely polarized}, if it is both polarized and term idempotent. The two conditions are independent. For example, the variety of groups is polarized but not term idempotent, and on the other hand, the variety $\mc{RS}$ of rectangular semigroups of Example~\ref{Ex:RS} is term idempotent but not polarized. However, if a variety satisfies both of these conditions, then their consequences are rather strong.

\begin{proposition}\label{P:ppident}
Let $\mcV$ be a purely polarized variety with a polar term $p(x)$. Let $u$ and $v$ be different terms of the type of $\mcV$. Then
\begin{equation*}
\mcV \models u = v\ \ \iff\ \ \mcV \models u = p(x) \ \ \mbox{and} \ \ \mcV \models v = p(x).
\end{equation*}
\end{proposition}
\begin{proof}
Suppose that $\mcV \models u = v$ for different terms $u$ and $v$. Since $\mcV$ is a term idempotent variety and the identity $u=v$ is nontrivial, both $u$ and $v$ are term idempotents of $\mcV$. Since $\mcV$ is polarized, they are $\mcV$-equivalent to $p(x)$. So $\mcV \models u = p(x)$ and $\mcV \models v = p(x)$.
\end{proof}

Typical examples of purely polarized varieties are the variety $\mc{C}$ of constant semigroups and, more generally, the varieties $\mc{C}_{\tau}$ of constant algebras of type $\tau$. Next we show a method of extending these examples.

\begin{example}
Constant semigroups may be defined by the consequences of the identities $x_1 \cdots x_n = y_1 \cdots y_n$ for all $n \geq 2$. Now, for $k \geq 2$, consider the variety~$\mc{C}_k$ of semigroups defined by the consequences of the identities $x_1 \cdots x_n = y_1 \cdots y_n$, where $n \geq k$. Note that $\mc{C}_2$ coincides with $\mc{C}$. The variety $\mc{C}_k$ satisfies all the identities $p = q$ of type $\tau$ such that neither $p$ nor $q$ is a product of less than $k$ variables. It is easy to see that both sides of every such identity are term idempotents of $\mc{C}_k$, and that $\mc{C}_k$ is polarized by the polar term $p(x) = x \cdots x$, where $x$ is repeated $k$ times. Hence $\mc{C}_k$ is purely polarized.

Recall the preordered set $(X\Om,\preceq)$ of Section~\ref{S:termidemp}. Consider the case when $\Om$ consists only of a single symbol of a binary operation. Note that the preorder $\preceq$ carries over to the free semigroup $X\mc{SG}$ over $X$. All term idempotents of $\mc{C}_k$ form an upper set of $(X\mc{SG},\preceq)$ that is generated by the element $x_1 \cdots x_k$.
\end{example}

Proposition~\ref{P:ppident} shows that both sides of every nontrivial identity satisfied in a purely polarized variety $\mcV$ are constant term idempotents of $\mcV$. Such identities will be called \emph{polar identities}.

\begin{corollary} \label{C:polaridentities}
A polarized variety $\mcV$ is purely polarized if and only if all nontrivial identities satisfied in $\mcV$ are polar.
\end{corollary}

A term $p(x)$ of type $\tau$ will be called a \emph{zero term} of a variety $\mcV$, if it is constant and for all $\omega \in \Om$ and every $1 \leq i \leq n$,
\begin{equation} \label{E:zeroomega}
\mcV \models\ \omega(x_1,\dots,x_{i-1},p(x),x_{i+1},\dots ,x_n) = p(x).
\end{equation}
Note that the identities of~\eqref{E:zeroomega} imply
\[
\mcV \models\ \omega(p(x),\dots,p(x)) = p(x).
\]
Hence a zero term is also a polar term, and any variety with a zero term is polarized. The pole of any algebra~$A$ in a variety with a zero term is the zero of $A$ (i.e. it forms a one-element sink of $A$).
If $p(x)$ is a zero term of $\mcV$, then the identities of~\eqref{E:zeroomega} may be easily generalized to all terms of type $\tau$. Indeed, for any term $t(x_1,\dots,x_n)$ and every $1 \leq i \leq n$,
\begin{equation}\label{E:zerot}
\mcV \models\ t(x_1,\dots,x_{i-1},p(x),x_{i+1},\dots,x_n) = p(x).
\end{equation}
If the type $\tau$ is plural, then there exists a term $t(x,y)$ involving both variables $x$ and $y$, such that the identities of~\eqref{E:zeroomega} imply
\[
\mcV \models\ p(x) = t(p(x),p(y)) = p(y).
\]
Thus if we restrict ourselves to plural types, a separate assumption that a zero term is constant is unnecessary.
The following proposition provides a basic property of polar identities.

\begin{proposition} \label{P:zeroterm}
Let $\mcV$ be a polarized variety. Let $u = v$ be a nontrivial polar identity true in $\mcV$. The following conditions are equivalent:
\begin{enumerate}
\item Every nontrivial consequence of a set of polar identities true in $\mcV$ is also polar.
\item Every nontrivial consequence of the identity $u = v$ is polar.
\item $\mcV$ has a zero term.
\end{enumerate}
\end{proposition}
\begin{proof}
$(1) \Rightarrow (2)$ This implication is obvious.

$(2) \Rightarrow (3)$ Assume (2). For any $\omega \in \Om$ and every $1 \leq i \leq n$, the identity
\begin{equation} \label{E:consequence1}
\omega(x_1,\dots,x_{i-1},u,x_{i+1},\dots ,x_n) = \omega(x_1,\dots,x_{i-1},v,x_{i+1},\dots ,x_n)
\end{equation}
is a nontrivial consequence of $u = v$. So, by (2), it is polar. The left-hand sides of the polar identities $u = v$ and \eqref{E:consequence1} are constant term idempotents. Hence they are $\mcV$-equivalent, because $\mcV$ is a polarized variety. Consequently, $\mcV$ satisfies the identity
\[
\omega(x_1,\dots,x_{i-1},u,x_{i+1},\dots ,x_n) = u.
\]
Therefore, the identities of \eqref{E:zeroomega} hold for the unary term $p(x) \coloneqq u(x,\dots,x)$, and so $p(x)$ is a zero term of $\mcV$.

$(3) \Rightarrow (1)$ Suppose that $\mcV$ has a zero term $p(x)$. (It is also a polar term of $\mcV$.) Consider a polar identity $u = v$ satisfied in $\mcV$. We will show that the consequences of this identity are also polar.

Let $t(x_1,\dots,x_n)$ be a term and let $1 \leq i \leq n$. Then the identity
\begin{equation} \label{E:consequence2}
t(x_1,\dots,x_{i-1},u,x_{i+1},\dots ,x_n) = t(x_1,\dots,x_{i-1},v,x_{i+1},\dots ,x_n)
\end{equation}
is a consequence of $u = v$. Since $u$ and $v$ are constant term idempotents, they are both $\mcV$-equivalent to $p(x)$. Thus, for the left-hand side of \eqref{E:consequence2}, we have
\[
t(x_1,\dots,x_{i-1},u,x_{i+1},\dots ,x_n) = t(x_1,\dots,x_{i-1},p(x),x_{i+1},\dots ,x_n) = p(x).
\]
An analogous identity holds for its right-hand side. Hence both sides of \eqref{E:consequence2} are constant term idempotents, and so \eqref{E:consequence2} is a polar identity.

If we substitute arbitrary terms for the variables of a constant term idempotent of $\mcV$, we again obtain a constant term idempotent. Thus an identity $u' = v'$, obtained from the polar identity $u = v$ by substituting some terms for its variables, is also polar.

If $u = v$ and $v = w$ are polar identities, then obviously $u = w$ is also a polar identity. Therefore all nontrivial consequences of any set of polar identities true in $\mcV$ are polar.
\end{proof}

\begin{corollary} \label{C:ppzero}
In a purely polarized variety, all polar terms are zero terms.
\end{corollary}
\begin{proof}
This follows directly by Proposition \ref{P:zeroterm} and Corollary \ref{C:polaridentities}.
\end{proof}

The following proposition provides an equational base for any purely polarized variety.

\begin{proposition}
A variety $\mcV$ of type $\tau$ is purely polarized if and only if it is defined by the identities
\begin{itemize}
\item[(a)] $p(x) = p(y)$,
\item[(b)] $\omega(x_1,\dots,x_{i-1},p(x),x_{i+1},\dots,x_n) = p(x)$ for all $\omega \in \Om$ and $1 \leq i \leq n$,
\item[(c)] $t_i = p(x)$ for all $i \in I$,
\end{itemize}
where $p(x)$ and $\{t_i \mid i\in I\}$ are some terms of type $\tau$.
\end{proposition}
\begin{proof}
$(\Rightarrow)$ Let $\mcV$ be a purely polarized variety with a polar term $p(x)$, defined by some identities $\{u_i = v_i \mid i \in I\}$. The polar term $p(x)$ satisfies the identity (a) and, by Corollary~\ref{C:ppzero}, also the identities of (b). By Proposition~\ref{P:ppident}, an identity $u_i = v_i$, for $i \in I$, has the same consequences as the identities $u_i = p(x)$ and $v_i = p(x)$. Thus we can construct an equational base of the required form.

$(\Leftarrow)$ Assume that $\mcV$ is defined by the identities of (a), (b) and (c). Then (a) and (b) imply that $p(x)$ is a polar term, and hence $\mcV$ is polarized. Thus, the defining identities (a), (b) and (c) of $\mcV$ are all polar identities. The identities of (b) show that $p(x)$ is a zero term of $\mcV$.
By Proposition~\ref{P:zeroterm}, all consequences of the defining identities are also polar, and so by Corollary~\ref{C:polaridentities}, $\mcV$ is purely polarized.
\end{proof}

We can see that purely polarized varieties form a very special class of algebras. Additionally, we will find that they interact with Mal'tsev products in an interesting way. First, let us look at what Theorem~\ref{T:newmain} says about the special case when $\mcW$ is a purely polarized variety. If we set $f(x,y,z)\coloneqq p(x)$ and $g(x,y,z)\coloneqq p(z)$ for a unary term $p(x)$, then the conditions (a), (b) and (c) reduce to: (a) $\mcV \models p(x) = x$, (b) $\mcW \models p(x) = p(y)$ and (c) $p(x)$ is a term idempotent of $\mcW$. We thus obtain the following corollary.

\begin{corollary} \label{C:5}
Let $\mcV$ and $\mcW$ be varieties of type $\tau$, and let $\mcW$ be purely polarized with a polar term $p(x)$. If $\mcV$ satisfies $p(x) = x$, then the Mal'tsev product $\mcV \circ \mcW$ is a variety.
\end{corollary}

\noindent Since every unary term different from a variable is a polar term of $\mc{C}_\tau$, additionally, we have the following.

\begin{corollary}
Let $\mcV$ be a variety of type $\tau$. Let $u(x)$ be a unary term of type $\tau$, which is not just a variable. If $\mcV$ satisfies $u(x) = x$, then the Mal'tsev product $\mcV \circ \mc{C}_\tau$ is a variety.
\end{corollary}

The assumption on the variety $\mcV$ is rather weak. For example, it is satisfied by all idempotent varieties and all strongly irregular varieties. Surprisingly, it is possible to prove a much more general result.
We conclude the paper by showing that the Mal'tsev product $\mcV \circ \mcW$ of any variety $\mcV$ and a purely polarized variety $\mcW$ is a variety.

This next lemma follows by Corollary~\ref{C:onesubalg} and the fact that each algebra in a purely polarized variety has exactly one idempotent.

\begin{lemma}\label{L:singlesubalg}
Let $\mcV$ and $\mcW$ be varieties of type $\tau$, and let $\mcW$ be a purely polarized variety. If $A \in \mcV \circ \mcW$, then exactly one  class of the $\mcW$-replica congruence of $A$ is a subalgebra of $A$, and all the other classes are singletons.
\end{lemma}

\begin{theorem}
Let $\mcV$ and $\mcW$ be varieties of type $\tau$. If $\mcW$ is a purely polarized variety, then the Mal'tsev product $\mcV \circ \mcW$ is a variety.
\end{theorem}
\begin{proof}
Let $A$ be a member of $\mcV \circ \mcW$. We will show that for any congruence~$\theta$ of $A$, the quotient algebra $A/\theta$ is also a member of $\mcV \circ \mcW$.

Let $\varrho$ be the $\mcW$-replica congruence of $A$.
By Lemma~\ref{L:singlesubalg}, exactly one class of $\varrho$, call it $E$, is a subalgebra of $A$. Moreover $E \in \mcV$, since $A \in \mcV \circ \mcW$. All the other classes of $\varrho$ are singletons.
Let $\bar{\theta} \coloneqq \varrho \vee \theta$, and let $a, b \in A$.
Since the class $E$ of $\varrho$ is a subalgebra of $A$ and other classes are singletons, it follows that
\begin{align*}\label{E:bartheta}
(a,b) \in \bar{\theta}\ \ \iff\ \  &\mbox{either}\ \ a/\theta \, \cap E \, \neq \varnothing \neq b/\theta\, \cap\, E\\
&\mbox{or}\ \ a/\theta = b/\theta \ \ \mbox{and}\ \ a/\theta\, \cap\, E = \varnothing.
\end{align*}
Hence, one class of $\bar{\theta}$, denoted by $E\theta$, contains $E$ and is a disjoint union of $\theta$-classes having a non-empty intersection with $E$, i.e.
\begin{equation*}
E\theta = \{a \in A \mid \exists e \in E,\ a \, \theta \, e\}
= \bigcup_{e \in E} e/\theta.
\end{equation*}
By the Second Isomorphism Theorem~\cite[Thm. 1.2.4]{RS02}, $E\theta$ is a subalgebra of~$A$ and
\begin{equation}\label{E:inV}
E\theta /(\theta \cap (E\theta)^2) \, \cong \, E/(\theta \cap E^2) \in \mcV.
\end{equation}
All the other classes of $\bar{\theta}$ coincide with $\theta$-classes which are disjoint from $E$.

Now recall that $(a,b) \in \bar{\theta}$ iff $(a/\theta, b/\theta) \in \bar{\theta}/\theta$, and
by the First Isomorphism Theorem~\cite[Thm. 1.2.3]{RS02},
\begin{equation}\label{E:FIT}
(A/\theta)/(\bar{\theta}/\theta) \cong A/\bar{\theta} \cong (A/\varrho)/(\bar{\theta}/\varrho).
\end{equation}
Since $A/\varrho$ belongs to $\mcW$, it follows that its quotient $(A/\varrho)/(\bar{\theta}/\varrho)$, and hence also the quotients $A/\bar{\theta}$ and $(A/\theta)/(\bar{\theta}/\theta)$, are members of $\mcW$.

Under the first isomorphism of~\eqref{E:FIT}, the $\bar{\theta}$-class $E\theta$ corresponds to the $\bar{\theta}/\theta$-class $\bar{E} \coloneqq \{e/\theta \mid e \in E\}$, the unique $\bar{\theta}/\theta$-class which is a subalgebra of $A/\theta$. Each of the remaining $\bar{\theta}/\theta$-classes consists of one $\theta$-class disjoint from~$E$. In particular, $\bar{\theta}/\theta$ is the $\mcW$-replica congruence of $A/\theta$. Note that, by~\eqref{E:inV}, $\bar{E} = E\theta /(\theta \cap (E\theta)^2) \in \mcV$. It follows that $A/\theta \in \mcV \circ \mcW$, and hence $\mf{H}(\mcV \circ \mcW) \subseteq \mcV \circ \mcW$. Therefore $\mcV \circ \mcW$ is a variety.
\end{proof}

\begin{corollary}
Let $\mcV$ be a variety of type $\tau$. Then the Mal'tsev product $\mcV \circ \mc{C}_\tau$ is a variety.
\end{corollary}

\end{document}